\documentclass{lms}

\usepackage{amsfonts,amssymb,amsmath}
\usepackage{graphicx}
\usepackage{xypic}
\usepackage[ps2pdf,                
bookmarks=true,
bookmarksnumbered=true,
hypertexnames=false,
breaklinks=true,
pdfborder={0 0 10}]{hyperref}

\def\Bbb{\mathbb}
\def\ca{{\mathcal A}}
\def\cc{{\mathcal C}}
\def\cd{{\mathcal D}}
\def\cf{{\mathcal F}}
\def\cg{{\mathcal G}}

\def\cs{{\mathcal S}}
\def\ct{{\mathcal T}}
\def\cu{{\mathcal U}}
\def\cx{{\mathcal X}}
\def\cy{{\mathcal Y}}

\def\cw{CW complex\ }
\def\gcw{$G$-CW complex\ }

\def\l{\label}
\def\|{|\bullet|}

\def\bz{\mathbb Z}
\def\bzc{\mathbb Z\cdot}

\def\bean{\begin{eqnarray}}
\def\eean{\end{eqnarray}}
\def\bea{\begin{eqnarray*}}
\def\eea{\end{eqnarray*}}
\def\nm{\nonumber}
\def\Bbb{\mathbb}

\def\eea{\end{eqnarray*}}

\def\os{\overset}
\def\er{\eqref}
\def\ubz{{\underline\bz}}
\def\mt{\mapsto}
\def\H{{\rm Hom}}

\def\a{\alpha}

\def\kp{\kappa}
\def\op{\oplus}
\def\ot{\otimes}
\def\gf{G\cf}
\def\ogf{\otimes_{\gf}}
\def\umt{{\underline M}_{tr}}

\newtheorem{theorem}{Theorem}[section] 
\newtheorem{lemma}[theorem]{Lemma}     
\newtheorem{corollary}[theorem]{Corollary}
\newtheorem{proposition}[theorem]{Proposition}

\newnumbered{assertion}{Assertion}    
\newnumbered{conjecture}{Conjecture}  
\newnumbered{definition}{Definition}
\newnumbered{hypothesis}{Hypothesis}
\newnumbered{remark}{Remark}
\newnumbered{note}{Note}
\newnumbered{observation}{Observation}
\newnumbered{problem}{Problem}
\newnumbered{question}{Question}
\newnumbered{algorithm}{Algorithm}
\newnumbered{example}{Example}
\newunnumbered{notation}{Notation} 

\title{A Functor Converting Equivariant Homology to Homotopy}

\author{Zhaohu Nie}

\classno{55N91}

\begin{document}
\maketitle

\begin{abstract}
In this paper, an equivariant version of the classical Dold-Thom theorem is proved. 
Let $G$ be a finite group, $X$ a $G$-space, and $k$ a covariant coefficient system on $G$. Then a topological abelian group ${\mathcal G}X\otimes_{G\cf} k$ is constructed by the coend construction. For a \gcw $X$, it is proved that there is a natural isomorphism $\pi_i({\mathcal G}X\otimes_{G\cf} k)\cong H_i^G(X;k)$, where the right hand side is the Bredon equivariant homology of $X$ with coefficients in $k$. At the end, several examples of this result are presented.
\end{abstract}


\section{Introduction}
From the point of view taken in this article, the classical Dold-Thom theorem concerns a functor which converts singular homology with $\Bbb Z$ coefficients to homotopy. This functor is the free abelian group functor, which sends a topological space $X$ to the free abelian group $\Bbb Z \cdot X$ generated by $X$ with a suitable topology. The Dold-Thom theorem \cite {dt} asserts that for a \cw  $X$,  
\bean
\pi_i(\Bbb Z \cdot X)\cong H_i(X;\Bbb Z).
\label{d-t}
\eean

More generally, given an abelian group $A$ and a topological space $X$, McCord \cite {mc} associated functorially to them a topological abelian group $B(X,A)$, which generalizes the construction of $\bzc X$. 
If $X$ is a CW complex, McCord \cite [Theorem 11.4]{mc} proved 
\bean
\pi_i(B(X,A))\cong H_i(X;A).
\l{anycoeff}
\eean

The goal of this paper is to give an equivariant version of \er{anycoeff}. 
Throughout the paper, we fix a finite group $G$. 
Let $X$ be a $G$-space. The equivariant analogue of the right hand side of \er{anycoeff} is the Bredon equivariant homology $H_i^G(X;k)$ of $X$ with coefficients in a covariant coefficient system $k$ on $G$. (We will recall the definitions of these in Section 2 (see also \cite b).) 
The analogue of $B(X,A)$ in \er{anycoeff} is our construction ${\mathcal G}X\otimes_{G\cf} k$ (see Definition \ref{xk} and cf. \er{analogy}), which functorially assigns a topological abelian group to a $G$-space $X$ and a covariant coefficient system $k$. We then prove in Section \ref{construction} our main theorem.
\begin{theorem}\l{main} For a $G$-CW complex $X$, one has a natural isomorphism
\bea
\pi_i({\mathcal G}X\otimes_{G\cf} k)\cong H_i^G(X;k),
\eea
where the right hand side is the Bredon equivariant homology of $X$ with coefficients in $k$.
\end{theorem}

The functor sending a $G$-space $X$ to the topological abelian group ${\mathcal G}X\otimes_{G\cf} k$ is the functor in the title, which converts equivariant homology to homotopy as shown in Theorem \ref{main}. 

The organization of this paper is as follows. In Section 2, we recall the definitions of covariant coefficient systems and Bredon equivariant homology. We also recall the coend construction here, which gives our construction in Section 3. 
We give the proof of Theorem \ref{main} in Section 4. Theorem \ref{main} generalizes certain results of Lima-Filho \cite{lf} and dos Santos \cite{ds} in the unstable case. We consider their relationship and some other examples in Section \ref{examples}. 


\section{Basic definitions}

In this section, we recall the definitions of covariant coefficient systems and Bredon equivariant homology. In the process, we also recall the definition of the coend construction. (See \cite{may} for more details.)


Let $G\cf$ be the category of finite $G$-sets and $G$-maps. By definition, a covariant coefficient system is a covariant functor $k:G\cf \to {\mathcal A}b$ to the category of abelian groups, which converts disjoint union to direct sum.

\begin{remark}\l{gfcg} Instead of $\gf$, one can consider the category $\mathcal G$ of orbit $G$-sets $G/H$, with $H$ a subgroup of $G$, and $G$-maps. Since each finite $G$-set $S$ is uniquely a disjoint union of orbits, and upon choosing a point, each orbit can be identified to $G/H$ for some $H$, it is clear that a covariant coefficient system can be equivalently defined to be a covariant functor $k:\cg\to \ca b$. Such equivalence happens throughout the paper (actually for all $\cg$-objects). Therefore the whole paper can be written in terms of $\cg$. However as we will see, working with $\gf$ gives us advantages, since it subsumes the category $\cf$ of finte sets, which is useful in view of Example \ref{e.g.}.
\end{remark}

We now recall the definition of Bredon equivariant homology following the singular approach of Illman \cite{i}, but we formulate it in the language of coends.
This formulation is suggestive to our later construction and proof.  

Recall that a simplicial object in a category $\mathcal C$ is a contravariant functor $\bigtriangleup^{op} \to \mathcal C$, where $\bigtriangleup$ is the category of sets $\underline n=\{0,1,\cdots,n\}$ and monotonic maps. For a topological space $X$, recall that the singular simplicial set of $X$ is defined as the represented functor
\bean\l{sx}
\cs X: \Delta^{op}\to \cs et;\ \underline n\mapsto \H_{\cu}(\triangle^n, X),
\eean
where $\cu$ is the category of topological spaces, and $\triangle^n$ is the standard topological $n$-simplex. 

Similarly, one defines a $\mathcal G$-object in a category $\mathcal C$, which we assume to have finite products, to be a contravariant functor $G\cf^{op} \to \mathcal C$, which converts disjoint union to product. For a $G$-space $X$, we define the associated $\cg$-space 
\bean\l{gx}
\cg X:G\cf^{op} \to \mathcal U;\ S\mapsto \H_{G\cu}(S,X),
\eean
as the represented functor in the category of $G$-spaces (where we regard $S$ as having the discrete topology). (In this paper, the notation 'Hom' for a topological category always denotes a morphism \emph{space}, unless a simplex $\Delta^n$ appears in the source.)


One can combine these two and define a $\cg$-simplicial object (or a simplicial $\cg$-object depending on context) in a category $\mathcal C$ to be a functor $G\cf^{op}\times \bigtriangleup^{op} \to \mathcal C$, which converts disjoint union in $G\cf$ to product. For a $G$-space $X$, we define the associated $\cg$-simplicial set as
\bean\l{gsx}
{\mathcal {GS}}X: G\cf^{op}\times \bigtriangleup^{op} \to {\mathcal S}et;\ S\times \underline n\mapsto \H_{G\cu}(S\times \bigtriangleup^n, X),
\eean
where $\bigtriangleup^n$ has the trivial $G$-action. 


The meaning of equivariancy (for homology) is best expressed in the coend construction which we now recall.

\begin{definition}\l{coend}
Let $\mathcal D$ be a small category, and $\mathcal C$ a category with finite products and 
all colimits. 
Let 
$$T:{\mathcal D}^{op} \to {\mathcal C}\ {\rm and}\ S: {\mathcal D} \to {\mathcal C}$$
be a contravariant and a covariant functor from ${\mathcal D}$ to ${\mathcal C}$, respectively. The coend of $T$ and $S$ 
is an object of $\cc$ which is the coequalizer of the following diagram
\bea
T\otimes _{\mathcal D} S={\rm coeq}\left(\coprod_{f:d\to e\in Mor(\cd)} T(e)\times S(d)\rightrightarrows \coprod_{d\in Ob(\cd)} T(d)\times S(d)\right),
\eea
where the two right arrows are $id\times f_*$ and $f^*\times id$ with $f_*=S(f)$ and $f^*=T(f)$. 
\end{definition}

Clearly, the coend construction is functorial: If $F:T\to T'$ and $G:S\to S'$ are two natural transformations, then one has a natural morphism
\bean\l{functoriality}
F\ot_\cd G: T\ot_\cd S\to T'\ot_\cd S'.
\eean


For the purpose of this paper, $\cd$ is either $\gf$
or $\cf$, the categories of finite $G$-sets or finite sets; $\cc$ is either $\cs et$ or $\cu$, the categories of sets or topological spaces. We also consider the natural forgetful functors $\ca b\to \cs et$ and $\ca b\to \cu$ with the discrete topology. 

When the objects of $\cc$ have elements (as for us), the coend $T\otimes_{\cd} S$ has the following explicit form
\bean\l{explicitcoend}
T\otimes_{\cd} S=\coprod_{d\in Ob(\cd)} T(d)\times S(d)/(\approx),
\eean
where the equivalence relation is generated by $tf^*\times s\approx t\times f_*s$ for a morphism $f:d\to e$ of $\cd$ and elements $t\in T(e)$ and $s\in S(d)$ whenever this makes sense. Here we write contravariant actions from the right to emphasize the analogy to tensor products. 

We give one example to illustrate the nature of the coend construction. 

\begin{example}\l{e.g.} Let $\cd=\cf$. 
A set $X$ gives rise to a contravariant functor (abusing notation)
$$X: {\cf}^{op} \to {\mathcal S}et;\ S\mapsto X^S=\H_{\cs et}(S,X).$$
An abelian group $A$ gives rise to a covariant functor 
\begin{gather}
A: \cf \to \ca b;\ S\mapsto A^S=\op_{s\in S} A,\nm\\
(f:S\to T)\mapsto \left(\op_{s\in S} A\to \op_{t\in T} A;\ (a_s)\mapsto (b_t=\sum_{s\in f^{-1}(t)} a_s)\right).\nm
\end{gather}

Then one has 
$$X\otimes _{\cf} A=B(X,A),$$
the abelian group generated by $X$ with coefficients in $A$. Actually in view of (\ref{explicitcoend}), 
$$X\otimes _{\cf} A=\coprod_{S} X^S\times A^S/(\approx).$$ 
If we define a map 
\bean\l{correspondence}
X^S\times A^S\ni ((x_s),(a_s))\mapsto \sum_{s\in S} a_s x_s\in B(X,A),
\eean
then the equivalence relation $(\approx)$ is exactly the one for identifications. 


Analogously one has corresponding constructions for simplicial sets 
and topological spaces. In particular, for a topological space $X$ one recovers $B(X,A)$ with its topology in \er{anycoeff} as
\bean\l{analogy}
B(X,A)=X\ot_\cf A=\coprod_{S} {\rm Hom}_\cu(S,X)\times A^S/(\approx),
\eean
where 
the product and coproduct are taken in the category $\cu$ of topological spaces.
\end{example}



Now returning to the definition of Bredon equivariant homology, we define the equivariant singular simplicial abelian group of $X$ with coefficients in $k$ as 
\bean\l{c_}
C^G_\bullet (X;k):={\mathcal {GS}}X \otimes_{\gf} k.
\eean
(It is intuitively clear that one gets a simplicial {\em abelian group} here. Also cf. Lemma \ref{abfun}.)

The $i$-th equivariant homology group of the $G$-space $X$ with coefficients in $k$ is defined to be the $i$-th homotopy group \cite{mays} of $C^G_\bullet(X;k)$, i.e.
\bean
H_i^G(X;k):=\pi_i(C^G_\bullet (X;k)).
\l{sdef}
\eean

There is a chain complex $C^G_* (X;k)$, associated to $C^G_\bullet (X;k)$, with the differentials as the alternating sums of the face maps. It is known \cite{mays} that the homotopy groups of $C^G_\bullet(X;k)$ are naturally isomorphic to the homology groups of $C^G_* (X;k)$. Therefore one has the following equivalent definition \cite{i}
\bea
H_i^G(X;k)=H_i(C^G_*(X;k)).
\eea

\section{Construction
}\l{hmtp}



In this section, we give our construction and list some simple properties. 

\begin{lemma}\l{abfun} Let $k$ be a covariant coefficient system. The coend construction gives a functor
$$\bullet\ot_{\gf} k: \cg\cu\to \ct\ca b,$$
from the category of $\cg$-spaces (as a functor category) to the category of topological abelian groups. 
\end{lemma}

\begin{proof} In view of \er{functoriality}, we only need to show the abelian group structure. 

For a $\cg$-space $\cx$, 
recall \er{explicitcoend} that 
$$\cx\ot_{\gf} k=\coprod_S \cx(S)\times k(S)/(\approx).$$
We denote the equivalence class of an element by $[-]$. 

We define the addition by juxtaposition: For $(x_S,\kappa_S)\in \cx(S)\times k(S)$ and $(x_T,\kappa_T)\in \cx(T)\times k(T)$, we define their sum to be $$[(x_S,\kappa_S)]+[(x_T,\kappa_T)]:=[((x_S,x_T),(\kappa_S,\kappa_T))]$$
where
$$((x_S,x_T),(\kp_S,\kp_T))\in (\cx(S)\times \cx(T))\times (k(S)\oplus k(T))=\cx(S\coprod T)\times k(S\coprod T).$$
We define the inverse by the inverse in $k$: For $(x,\kp)\in \cx(S)\times k(S)$, we define 
$$-[(x,\kp)]:=[(x,-\kp)].$$
Now let's check compatibility. By our definition, 
$$[(x,\kp)]+[(x,-\kp)]=[((x,x),(\kp,-\kp))]$$
for $((x,x),(\kp,-\kp))\in \cx(S\coprod S)\times k(S\coprod S).$
Consider the folding map 
$$\nabla=id\coprod id: S\coprod S\to S.$$ 
One sees that 
$$((x,x),(\kp,-\kp))=(x\nabla^*,(\kp,-\kp))\approx (x,\nabla_*(\kp,-\kp))\in \cx(S)\times k(S).$$
It is clear 
that $\nabla_*(\kp,-\kp)=\kp+(-\kp)=0$. It is also clear that $[(x,0)]=0$ in $\cx\ot_{\gf} k$ (since 0 ``comes'' from $\emptyset$). 

For a natural transformation $\phi:\cx\to \cy$, let's check that \er{functoriality}
$$\phi_*:=\phi\ogf id:\cx\ogf k\to \cy\ogf k$$ 
is a homomorphism. In view of our definition of the addition, this boils down to the commutativity of following diagram
$$\xymatrix{
\cx(S\coprod T)\ar[r]^=\ar[d]^{\phi(S\coprod T)} & \cx(S)\times \cx(T)\ar[d]^{\phi(S)\times \phi(T)}\\
\cy(S\coprod T)\ar[r]^= & \cy(S)\times \cy(T).
}$$
\end{proof}

\begin{definition}\l{xk} For a $G$-space $X$ and a covariant coefficient system $k$, we define the topological abelian group $\cg X\ot_{G\cf} k$ as the coend construction in view of \er{gx} and Lemma \ref{abfun}. 
\end{definition}

Now we prove a simple lemma which will be used later. 

A $\mathcal G$-homotopy from $\cg$-spaces ${\mathcal X}$ to ${\mathcal Y}$ is a $\mathcal G$-map (natural transformation)\\
$H: I\times {\mathcal X} \to {\mathcal Y}$, where $I$ is the unit interval, and $I\times {\mathcal X}$ is the $\mathcal G$-space defined by $(I\times {\mathcal X})(S)=I\times {\mathcal X}(S)$ and similarly for morphisms. Two $\mathcal G$-maps $\phi,\psi: {\mathcal X}\to {\mathcal Y}$ are said to be $\mathcal G$-homotopic if there exists a $\mathcal G$-homotopy $H$ from ${\mathcal X}$ to ${\mathcal Y}$ such that $H(0, \cdot)=\phi$ and $H(1, \cdot)=\psi$. A $\mathcal G$-map $\phi: {\mathcal X}\to {\mathcal Y}$ is called a $\mathcal G$-homotopy equivalence if there exist a $\mathcal G$-map $\phi':{\mathcal Y}\to {\mathcal X}$ such that $\phi'\circ \phi$ is $\mathcal G$-homotopic to $id_{\mathcal X}$ and $\phi\circ \phi'$ is $\mathcal G$-homotopic to $id_{\mathcal Y}$. We have the following lemma. 

\begin{lemma}\l{prehom} If $\phi, \psi: {\mathcal X}\to {\mathcal Y}$ are $\mathcal G$-homotopic, then $\phi_*, \psi_*:{\mathcal X}\otimes_{G\cf} k \to {\mathcal Y} \otimes_{G\cf} k$ are homotopic through homomorphisms.
\end{lemma}

\begin{proof} Suppose that $H: I\times {\mathcal X}\to {\mathcal Y}$ is a $\mathcal G$-homotopy such that $H(0,\cdot)=\phi$ and $H(1,\cdot)=\psi$. Then by Lemma \ref{abfun}, $H$ induces a continuous homomorphism 
$$H_*: (I\times {\mathcal X})\otimes_{G\cf} k \to {\mathcal Y} \otimes_{G\cf} k.$$

We also have a natural map 
$$i: I\times ({\mathcal X}\otimes_{G\cf} k)\to (I\times {\mathcal X})\otimes_{G\cf} k;\ (t, [(x, \kp)])\mt [((t, x),\kp)],$$
where $(x,\kp)\in \cx(S)\times k(S)$ for some $S$. 

Denote the composition of $H_*$ and $i$ by ${\mathcal H}: I\times ({\mathcal X}\otimes_{G\cf} k)\to {\mathcal Y} \otimes_{G\cf} k$. Then $ {\mathcal H}(0, \cdot)=\phi_*$, and ${\mathcal H}(1, \cdot)=\psi_*$. 
\end{proof}

From Lemmas \ref{prehom} and \ref{abfun}, we have the following corollary.

\begin{corollary}\l{homotopy} If $\phi: {\mathcal X}\to {\mathcal Y}$ is a $\mathcal G$-homotopy equivalence, then $\phi_*: {\mathcal X}\otimes_{G\cf} k \to {\mathcal Y} \otimes_{G\cf} k$ is a homotopy equivalence. 
\end{corollary}

\section{Proof of the main theorem}\l{construction}

We first prove three lemmas concerning the geometric realization $\|$ (see \cite{mays} for reference), and then use them to prove our Theorem \ref{main}. 

In view of \er{gsx}, we know that $|\cg\cs X|$ is a $\cg$-CW complex, whose value on a finite $G$-set $S$ is
$$|\cg\cs X|(S)=|\cg\cs X(S,\bullet)|=|(\underline n\mapsto \H_{G\cu}(S\times \Delta^n,X))|.$$

\begin{lemma} \l{1} Let $\cg\cs X$ be the $\cg$-simplicial set associated to a $G$-space $X$. We have the following homeomorphism 
$$|{\mathcal {GS}}X\otimes_{G\cf}k|\simeq |{\mathcal {GS}}X|\otimes_{G\cf}k,$$
i.e. the geometric realization commutes with the coend construction. 
\end{lemma}

\begin{proof} We proceed by the following sequence of homeomorphisms: 
\begin{align}
 &|{\mathcal {GS}}X\otimes_{G\cf}k|
\nm\\
=&\left|{\rm coeq}\left(\coprod \cg\cs X(S,\bullet)\times k(T)\rightrightarrows \coprod \cg\cs X(S,\bullet)\times k(S)\right)\right| 
\nm\\
\simeq&{\rm coeq}\left(\coprod |\cg\cs X(S,\bullet)|\times k(T))\rightrightarrows \coprod |\cg\cs X(S,\bullet)|\times k(S)\right) 
\nm\\
=&|{\mathcal {GS}}X|\otimes_{G\cf}k, 
\nm
\end{align}
where the first and the last equalities follow from Definition \ref{coend}, and the second from the following. 

It is well known \cite{mays} that $\|$ is the left adjoint of the singular simplicial set functor $\cs$ in \er{sx}. Therefore $\|$ commutes with all colimits. The second equality follows from this since $\times k(S)$ (with the discrete topology), coproducts and the coequalizer are all colimits. 
\end{proof}

For a $G$-space $X$, $\cs X$ \er{sx} is a simplicial $G$-set, and therefore $|\cs X|$ is a $G$-CW complex. The associated $\cg|\cs X|$ \er{gx} is a $\cg$-CW complex, whose value on a finite $G$-set $S$ is 
$$\cg|\cs X|(S)=\H_{G\cu} (S,|\cs X|).$$

\begin{lemma} \l{1.5} For a $G$-space $X$, one has an isomorphism of $\cg$-CW complexes
$$|\cg\cs X|\simeq \cg |\cs X|.$$
\end{lemma}

\begin{proof} In this proof, we switch our language to the orbit category $\cg$ in view of Remark \ref{gfcg}, since the geometric realization commutes with finite products \cite{mays} (we work in the category of compactly generated spaces).


Fix an orbit $G/H$ and observe \er{gsx} that 
\begin{align}
\cg\cs X(G/H,\underline n)=\H_{G\cu}(G/H\times \bigtriangleup^n,X)=\H_{\cu}(\bigtriangleup^n,X^H)=\cs X^H(\underline n),\nm
\end{align}
where $X^H$ is the fixed point set of $X$ by $H$. Therefore as a simplicial set,
\bea
{\mathcal {GS}}X(G/H,\bullet)={\mathcal S}X^H.
\eea

It is clear that $\cs X^H=(\cs X)^H$, the fixed point set by $H$ of the simplicial $G$-set $\cs X$. One has $|(\cs X)^H|=|\cs X|^H$ by \cite[eqn. (V.1.3)]{may}.

Therefore 
\begin{align}
|\cg\cs X|(G/H)=|\cg\cs X(G/H,\bullet)|=|\cs X^H|=|\cs X|^H=\cg |\cs X|(G/H).\nm
\end{align}
\end{proof}

\begin{lemma} \l{2} Let $X$ be a $G$-CW complex. Then there is a natural $\cg$-homotopy equivalence
$$\cg j: \cg|\cs X|\to \cg X.$$
\end{lemma}

\begin{proof} One has a natural map of $G$-spaces \cite{mays}
\bea
j:|\cs X|\to X,
\eea
which is a $G$-weak equivalence, i.e. $j^H:|\cs X|^H=|\cs X^H|\to X^H$ is a weak equivalence for any $H$. 

$|\cs X|$ is a $G$-CW complex. If $X$ is a $G$-CW complex, then by the equivariant Whitehead theorem \cite[Cor. I.3.3]{may}, $j$ is a $G$-homotopy equivalence. 


Now by the functoriality of \er{gx}, we see that 
$$\cg j: \cg |\cs X|\to \cg X$$
is a $\cg$-homotopy equivalence. 
\end{proof}

Now we are ready for the proof of Theorem \ref{main}. 

\begin{proof}[\thinspace of Theorem \ref{main}] First recall that for a Kan simplicial set $K_\bullet$, e.g. the simplicial abelian group $\cg\cs X\otimes_{G\cf} k$, there is a natural isomorphism \cite{mays}
\bean\l{referee}
\pi_i(K_\bullet)\cong \pi_i(|K_\bullet|).
\eean

We proceed by the following sequence of natural isomorphisms:
\begin{align}
H_i^G(X;k)
          &=\pi_i(C^G_\bullet (X;k)) && {\rm(definiton\ \er{sdef})}\nm\\
          &=\pi_i({\mathcal {GS}}X\otimes_{\gf}k) && {\rm(definition\ \er{c_})}\nm\\
          &\cong\pi_i(|{\mathcal {GS}}X\otimes_{\gf} k|) && {\rm(equation\ \er{referee})}\nm\\
          &\cong\pi_i(|{\mathcal {GS}}X|\otimes_{\gf} k) && {\rm(Lemma\ \er 1)}\nm\\
          &\cong\pi_i(\cg|{\mathcal {S}}X|\otimes_{\gf} k) && {\rm(Lemma\ \er {1.5})}\nm\\
          &\cong\pi_i({\mathcal G}X\otimes_{\gf} k).&& {\rm(Lemma\ \er 2 \ and \ Corollary \ \er{homotopy})}\nm
\end{align}
\end{proof}



\begin{remark} In a joint work \cite{dsn} with P. dos Santos, we will investigate what happens when the covariant coefficient system comes from a Mackey functor (see \cite{may} for the definition). It turns out that our construction in Definition \ref{xk} will then have more structure. We then apply our construction and the knowledge from \cite{ds} to study the Eilenberg-MacLane spectrum and $RO(G)$-graded homology associated to the Mackey functor (again see \cite{may} for the definition). 
\end{remark}

\section{Examples}\l{examples}

In this section, we calculate several examples of our construction $\cg X\ogf k$ for simple covariant coefficient systems, and then compare our Theorem \ref{main} to some other results. 

Let $\underline \bz$ be the constant (on orbits) covariant coefficient system at $\bz$, i.e. 
$$\underline \bz: \cg\to \ca b;\ G/H\mapsto \bz,\ (f:G/H\to G/K)\mapsto (id:\bz\to \bz),$$
and it converts disjoint union to direct sum. 

More categorically, letting 
$$q:\gf\to \cf;\ S\mapsto S/G$$
be the quotient functor, one actually has
$$\underline\bz=\bz\circ q: G\cf\to \cf\to \ca b,$$
where $\bz$ is as in Example \ref{e.g.}. 

\begin{proposition}\l{prop:6.1} For a $G$-space $X$, one has a natural homeomorphism 
$$\cg X\otimes_{G\cf} \underline \bz\simeq \bzc (X/G).$$
\end{proposition}

\begin{proof} We first define natural maps in both directions. 

Let $q$ (abusing notation) be the following composition (cf. \er{analogy}):
\begin{align}
q:&\cg X\otimes_{\gf} \underline \bz=\coprod_{S\in \gf} \H_{G\cu}(S,X)\times \bz^{q(S)}/(\approx)\nm\\
\os{q\times id}\longrightarrow & \coprod_{q(S)\in \cf} \H_\cu(q(S),X/G)\times \bz^{q(S)}/(\approx)\to (X/G)\ot_\cf \bz=\bz\cdot (X/G).\nm
\end{align}

For the other direction, let $\a\in \H_\cu(T,X/G)$ with $T\in \cf$ and consider the following pullback diagram in the category $G\cu$ of $G$-spaces ($T$ and $X/G$ have trivial $G$-actions):
$$\xymatrix{
p(T)\ar[r]^{p(\a)}\ar[d] & X\ar[d]^\pi\\
T\ar[r]^\a & X/G,
}$$
where $\pi:X\to X/G$ is the quotient map. One sees that $p(T)\in \gf$ since $G$ is finite. 
Note that $q(p(T))=T$ (since $q(\pi)=id$). 
Now consider the following composition 
\begin{align}
p:&\bz\cdot (X/G)=(X/G)\ot_\cf \bz=\coprod_{T\in \cf} \H_\cu(T,X/G)\times \bz^T/(\approx)\nm\\
\os{p\times id}\longrightarrow & \coprod_{p(T)\in\gf} \H_{G\cu}(p(T),X)\times \bz^{q(p(T))}/(\approx)\to \cg X\ogf \underline \bz.\nm
\end{align}

One can check that both $q$ and $p$ are well defined, and they are inverses of each other: Clearly $q\circ p=id$; the natural map $S\to p(q(S))$, by the universality of $p(q(S))$, and the equivalence relation give $p\circ q=id$. 
\end{proof}



For a $G$-CW complex $X$, our Theorem \ref{main}, Proposition \ref{prop:6.1} and the Dold-Thom theorem \er{d-t} together give
$$
H_i^G(X,\ubz)\cong  \pi_i(\cg X\otimes_{\gf}\underline \bz)\cong \pi_i(\bz\cdot (X/G))\cong     H_i(X/G,\bz),
$$
which recovers an easy and well-known fact
\cite[page 35]{may}. 

For a (left) $G$-module $M$, i.e. an abelian group $M$ with a left additive $G$-action, $M\otimes X=B(X,M)$ has a naturally induced $G$-action by 
$$g(\sum m_i x_i)=\sum (gm_i) (gx_i).$$ 
Dos Santos \cite {ds} proved that when $X$ is a $G$-CW complex, 
$$\pi_i((M\otimes X)^G)=H_i^G(X;\underline {M}_{tr}),$$
where ${\underline M}_{tr}$ is the covariant coefficient system defined by 
\begin{gather}
\umt:\gf\to \ca b;\ S\mapsto \H_{G}(S,M),\nm\\
(f:S\to T)\mapsto \left(\H_{G}(S,M)\to \H_{G}(T,M);\ \psi\mapsto (t\mapsto \sum_{s\in f^{-1}(t)} \psi(s))\right).\nm
\end{gather}
Before this, Lima-Filho \cite{lf} studied the special case when $M=\bz$ with a trivial $G$-action. 
Our Theorem \ref{main} covers this result in view of the following.

\begin{proposition}  
For a $G$-space $X$, one has a natural homeomorphism 
$$\cg X\otimes_{\gf} {\underline M}_{tr}\simeq (M\otimes X)^G.$$
\end{proposition}

\begin{proof} In view of \er{explicitcoend}, one has 
\begin{align}
\cg X\otimes_{\gf} {\underline M}_{tr}=&\coprod_{S\in \gf} \H_{G\cu}(S,X)\times \H_{G}(S,M)/(\approx)\nm\\
=&\coprod_{S\in \gf} \H_{G\cu}(S,X\times M)/(\approx).\nm
\end{align}
In view of Example \ref{e.g.}, one has
\begin{align}
M\otimes X=B(X,M)=&\coprod_{T\in \cf} \H_\cu(T,X)\times \H(T,M)/(\approx)\nm\\
=&\coprod_{T\in \cf} \H_\cu(T,X\times M)/(\approx).\nm
\end{align}

Forgetting the $G$-actions, one has a natural map $\cg X\otimes_{\gf} \umt\to M\otimes X$, which factorizes through
$$f:\cg X\otimes_{\gf} \umt\to (M\otimes X)^G,$$
since in view of \er{correspondence} $f([(\alpha,\beta)])=\sum_{s\in S} \beta(s)\alpha(s)$ is invariant under $G$, for $(\alpha,\beta)\in \H_{G\cu}(S,X)\times \H_{G}(S,M)$. 

Now we want to show that there is a natural map 
$$h:(M\otimes X)^G\to \cg X\otimes_{\gf} \umt.$$




For an element $a$ in $M\ot X$, one can always choose a representative 
\bean\l{rep}
(\gamma,\delta)\in \H_\cu(T,X)\times \H(T,M)=\H_\cu(T,X\times M)
\eean 
for some $T$ with $\gamma$ injective. (If not, use ${\rm Im}\gamma$ instead and apply the equivalence relation. This amounts to adding coefficients of similar terms in view of \er{correspondence}. Clearly the same applies to $\cg X\ogf \umt$.) 

In particular if $a\in (M\ot X)^G$, we see that ${\rm Im}(\gamma,\delta)\subset X\times M$ is invariant under the $G$-action and thus a finite $G$ set with the induced $G$-action. Since $(\gamma,\delta):T\to {\rm Im}(\gamma,\delta)$ is a bijection by assumption, one sees that $T$ has a natural $G$-action, such that 
$$(\gamma,\delta)\in \H_{G\cu}(T,X\times M).$$ 
We then define
$$h(a)=[(\gamma,\delta)].$$

It is clear that $f$ and $h$ are well defined and inverses of each other (upon choosing representatives as in \er{rep}). 
\end{proof}

\begin{acknowledgements}
The author would like to express his gratitude to Professor Blaine Lawson for his interest, discussions and encouragements. He thanks Christian Haesemeyer for stimulating discussions and a careful reading of an early version of the manuscript. He also thanks Pedro dos Santos and Paulo Lima-Filho for useful discussions, and the referee for his/her nice comments. 
\end{acknowledgements}

\affiliationone{
   Zhaohu Nie\\
   Department of Mathematics\\ 
   Texas A\&M University\\ 
   College Station\\ TX 77843-3368\\ 
   USA
   \email{nie@math.tamu.edu}}

\end{document}